\documentclass[12pt,a4paper]{article}
\usepackage{amsmath,amssymb,amsfonts}
\def\Bbb{\mathbb}

\title{\bf Linear relations among asymptotic frequencies in continued fractions}

\author{Kurt Girstmair}
\date{}

\makeatletter
\let\@@maketitle=\maketitle
\def\maketitle{\def\thispagestyle##1{\relax}\@@maketitle}
\makeatother
\newtheorem{theorem}{Theorem}
\newtheorem{prop}{Proposition}
\newtheorem{lemma}{Lemma}

\def\BE{\begin{equation}}
\def\EE{\end{equation}}
\def\BD{\begin{displaymath}}
\def\ED{\end{displaymath}}
\def\BA{\begin{array}}
\def\EA{\end{array}}
\def\BEA{\begin{eqnarray*}}
\def\EEA{\end{eqnarray*}}
\def\BI{\bibitem}

\def\N{\Bbb N}
\def\Z{\Bbb Z}
\def\Q{\Bbb Q}
\def\R{\Bbb R}

\def\H{\Bbb H}

\def\phi{\varphi}

\def\MB{\mbox}
\def\LD{\ldots}

\def\sminus{\smallsetminus}
\def\DIV{\,|\,}

\def\MN{\medskip\noindent}
\def\STOP{\hfill$\Box$}

\def\H{{\cal H}}
\def\M{{\cal M}}

\begin{document}
\maketitle

\begin{abstract}
\noindent
Let $H(m,d)$ denote the asymptotic frequency of the natural numbers $k\equiv d \mod m$ in the continued fraction expansions
of almost all numbers $x\in[0,1)$. For a fixed number $m\ge 4$, we study $\Q$-linear relations among the numbers $H(m,d)$, $1\le d\le m-3$,
i.e., vectors $(c_1,\LD,c_{m-3})\in\Q^{m-3}$ such that
\BD
  \sum_{d=1}^{m-3} c_dH(m,d)=0.
\ED
We restrict ourselves to the symmetric case $c_d=c_{m-2-d}$. In the end, we obtain a basis of the $\Q$-vector space of these relations for prime powers $m$
and for $m=pq$, where $p\ne q$ are primes.
\end{abstract}

\section*{1. Introduction}
Let $x$ be a real number, $0\le x<1$, and $x=[0,a_1,a_2,a_3,\LD]$ the (regular) continued fraction expansion of $x$.
Birkhoff's ergodic theorem (see \cite[p. 221]{IoKr}) allows the computation of the asymptotic frequency of a set $\M$ of natural numbers
among the digits $a_j$  of $x$ for almost all $x\in[0,1)$. For instance, if $\M=\{m\}$, then
\BD
\label{1.0}
   \lim_{n\to\infty} \frac{|\{j;\,j\le n, a_j=m \}|}n=\log_2\left(\frac{{(m+1)}^2}{m(m+2)}\right)
\ED
for all $x\in[0,1)$ up to those in a set of Lebesgue measure $0$. Here $\log_2$ means the logarithm to the base $2$.

In this paper we deal with $\M=\{k\in\N;\, k\equiv d\mod m\}$, where $m$ and $d$ are natural numbers, $m\ge 2$, $d\le m$.
In this case a meanwhile classical result says
\BE
\label{1.1}
  \lim_{n\to\infty} \frac{|\{j;\,j\le n, a_j\equiv d\mod m \}|}n=\log_2\left(\frac{\Gamma(d/m)\Gamma((d+2)/m)}{\Gamma((d+1)/m)^2}\right),
\EE
for all $x\in[0,1)$ up to those in a set of Lebesgue measure $0$ (see \cite{No} and \cite[p. 227]{IoKr}). Here $\Gamma(\LD)$ denotes Euler's gamma function.
We write $H(m,d)$ for the number in (\ref{1.1}).

One has $H(4,1)=1/2$ (see \cite[p. 228]{IoKr}). In the remaining cases the numbers $H(m,d)$ are supposed to be transcendental (a partial result in this direction is \cite[Satz 8]{Gi}).
However, the arithmetic properties of these numbers are still awaiting a more thorough investigation. So far our previous paper \cite{Gi} seems to be the only attempt of this kind.
In the said article we exhibited certain inhomogeneous $\Q$-linear relations among the numbers $H(m,d)$, for instance,
\BD
  H(12,1)=H(12,6)+H(12,7)+H(12,8)+1/3.
\ED
Here we study, in a more systematic way, homogeneous $\Q$-linear relations, i. e.,
rational numbers $c_d$ such that
\BD
\label{1.3}
  \sum_{d=1}^m c_d\,H(m,d)=0.
\ED

We fix $m$ for the time being and simply write $H_d$ for $H(m,d)$, $1\le d\le m$.
In this paper we investigate the aforesaid numbers $c_d$ under the following restrictions: First, $c_m=c_{m-1}=c_{m-2}=0$, and, second, $c_d=c_{m-2-d}$ for $1\le d\le m-3$. The latter condition makes
$(c_1, c_2, \LD, c_{m-3})$ {\em symmetric} with respect to $m/2-1$. Of course, it amounts to the same if we introduce
\BD
\label{1.5}
 S_d=H_d+H_{m-2-d}, \enspace 1\le d< m/2-1,
 \ED
 together with
\BD
\label{1.7}
 S_{m/2-1}=H_{m/2-1},
\ED
if m is even,
and study the possible numbers $c_d\in \Q$, $1\le d\le m/2-1$ , that yield
\BD
\label{1.9}
   \sum_{1\le d\le m/2-1} c_d\,S_d=0.
\ED
We call the numbers $S_d$, $1\le d\le m/2-1$, the {\em symmetric frequencies modulo} $m$.
The symmetric frequencies are logarithms of numbers in a cyclotomic field. This allows the application of results about cyclotomic fields (in particular, about units) to our problem.

Put $m'=\lfloor m/2\rfloor$. By means of the formula $\Gamma(z)\Gamma(1-z)=\pi/\sin(\pi z)$, $0<z<1$, we obtain
\BE
\label{1.11}
  S_d=\log_2\left(\frac{\sin(\pi(d+1)/m)^2}{\sin(\pi d/m)\sin(\pi(d+2)/m)}\right),\enspace d=1,\LD, m'-2,
\EE
and
\BE
\label{1.15}
 S_{m'-1}=\log_2\left(\frac {\sin (\pi m'/m)}{\sin(\pi(m'-1)/m)}\right)
\EE
(see \cite[formulas (18), (26), (31)]{Gi}). Observe that all logarithms are taken from positive real numbers. These formulas make the connection of symmetric frequencies
with cyclotomic fields obvious.

In order to study linear relations among the numbers $S_d$, we need some framework.
We work with the polynomial ring $\Q[X_1,\LD, X_{m'-1}]$ in $m'-1$ variables.
An $S$-{\em relation $R$} is a linear form $R=\sum_{k=1}^{m'-1}c_kX_k$, $c_k\in\Q$, such that $R(S_1,\LD,S_{m'-1})=0$.
Of course, the $S$-relations form a $\Q$-linear subspace of the space
$\Q X_1+\LD +\Q X_{m'-1}$ of linear forms. Our main goal is the description of a $\Q$-basis of this subspace.

To this end we introduce numbers $U_k, k=1,\LD, m'$, which are closely connected with the $S_d$'s but more appropriate for this purpose.
Let
\BE
\label{1.17}
 U_k=\log_2\left(\frac{\sin(\pi k/m)}{\sin(\pi/m)}\right) \MB{ for } k=1,\LD, m'.
\EE
Note that $U_1=0$.
In Section 2 we will see that the $\Q$-vector spaces $\Q S_1 +\LD +\Q S_{m'-1}$ and $\Q U_2+\LD+\Q U_{m'}$ are the same. More precisely, we express the generators $(S_1,\LD,S_{m'-1})$ and $(U_2,\LD, U_{m'})$
as linear combinations in terms of each other. This allows reducing $S$-relations to $\Q$-linear relations among the $U_k$'s.

For this purpose we need the space  $\Q Y_2+\LD+\Q Y_{m'}$ of linear forms in the variables $Y_2,\LD,Y_{m'}$. As in the case of the $S_d$'s, a $U$-{\em relation}
is a linear form $R'=\sum_{k=2}^{m'}c_k'Y_k$
such that $R'(U_2,\LD,U_{m'})=0$.
The results of Section 2 yield a linear isomorphism
\BD
   \phi: \Q Y_2+\LD+\Q Y_{m'} \to \Q X_1+\LD+\Q X_{m'-1}
\ED
which maps the $\Q$-subspace of $U$-relations onto the subspace of $S$-relations.

In Section 3 we describe a $\Q$-basis of the space of $U$-relations in the case of a prime power $m=p^n$. By means of $\phi$ we know a $\Q$-basis of the space of $S$-relations
for these numbers $m$ in every particular case.

The case of composite numbers $m$ seems to be quite involved. We settle only the cases $m=pq$, where $p$ and $q$ are two distinct primes $\ge 3$ (Section 4),
and $m=2p$, $p$ a prime $\ge 3$ (Section 5).

In Section 6 we consider a family of short $S$-relations. Thereby, we show that in many cases not every subset of $S_1,\LD,S_{m'-1}$ of appropriate size can serve as a basis of the space $\Q S_1+\LD+\Q S_{m'-1}$.

\section*{2. The connection between the numbers $S_d$ and $U_k$}

Let $m'=\lfloor m/2\rfloor$ be as above. Obviously, we have, from (\ref{1.11}), (\ref{1.15}), and (\ref{1.17}),
\BE
\label{2.1}
S_d=2U_{d+1}-U_d-U_{d+2}\enspace \MB{ for } \enspace d=1,\LD,m'-2,
\EE
and
\BE
\label{2.3}
S_{m'-1}=U_{m'}-U_{m'-1}.
\EE
Conversely, we have
\BE
\label{2.5}
 U_{k+1}=\sum_{d=1}^k\sum_{j=d}^{m'-1} S_j, \MB{ for } k=1,\LD, m'-1.
\EE
This identity is a matter of telescoping sums. Indeed,
\BD
 \sum_{j=d}^{m'-2}(2U_{j+1}-U_j-U_{j+2})+U_{m'}-U_{m'-1} =U_{d+1}-U_d
\ED
for $d=1,\LD, m'-1$. Now (\ref{2.5}) follows from
\BD
 \sum_{d=1}^k(U_{d+1}-U_d)=U_{k+1}-U_1=U_{k+1},
\ED
$k=1,\LD,m'-1$.
The identity (\ref{2.5}) can also be written
\BE
\label{2.7}
  U_{k+1}=\sum_{d=1}^k d\, S_d +\sum_{d=k+1}^{m'-1}k S_d.
\EE
In order to apply the identities (\ref{2.1}), (\ref{2.3}), and (\ref{2.7}) to linear relations,
we define a $\Q$-linear isomorphism
\BD
\label{2.9}
  \phi:\Q Y_2+\LD+\Q Y_{m'} \to \Q X_1+\LD+\Q X_{m'-1}: Y_{k+1}\mapsto \sum_{d=1}^k d\, X_d +\sum_{d=k+1}^{m'-1}k X_d,
\ED
$k=1,\LD,m'-1$. Of course, the inverse mapping of $\phi$ is given by
\BD
 X_d\mapsto\begin{cases} 2Y_{d+1}-Y_d-Y_{d+2}, & d=1,\LD, m'-2;\\
                         Y_{m'}-Y_{m'-1},       & d=m'-1.
           \end{cases}
\ED
Here and in the sequel we use the convention $Y_1=0$.
From the above it is obvious that $\phi$ maps a $U$-relation $R'=\sum_{k=2}^{m'}c_k'Y_k$ onto an $S$-relation $R=\sum_{k=1}^{m'-1}c_kX_k$.
Indeed, $\phi$ induces an isomorphism between the two spaces of $\Q$-linear relations. Since we are going to describe a basis of the space of
$U$-relations, it is advisable to know how $\phi$ transforms a $U$-relation $R'$ as above. We obtain
\BE
\label{2.11}
  \phi\left(\sum_{k=2}^{m'}c_k'Y_k\right)=\sum_{d=1}^{m'-1} c_d X_d\enspace  \MB{ with }\enspace c_d=\sum_{k=2}^d(k-1)c_k'+d\sum_{k=d+1}^{m'} c_k'.
\EE
Note that the coefficients $c_d$ can be computed, in a convenient way, by recursion. Indeed, we have $c_d=e_d+df_d$ with $e_1=0$, $f_1=\sum_{k=2}^{m'}c_k'$ and
$e_{d+1}=e_d+dc_{d+1}'$, $f_{d+1}=f_d-c_{d+1}'$, $d=1,\LD,m'-2$.

\section*{3. The case $m=p^n$}

As above, let $m\ge 4$ and $m'=\lfloor m/2\rfloor$.
Let $\zeta_m=e^{2 \pi i/m}$ be the standard primitive $m$th-root of unity. Then
\BE
\label{3.1}
  \frac{\sin(\pi k/m)}{\sin(\pi/m)} =\zeta_{2m}^{1-k}\frac{1-\zeta_m^k}{1-\zeta_m}
\EE
for every $k\in \Z$.

Let $m=p^n$ be a prime power.
We consider the multiplicative subgroup $M$ of $\R^{\times}$ that is generated by the positive numbers
\BD
  \frac{\sin(\pi k/m)}{\sin(\pi/m)},\enspace k=2,\LD, m'.
\ED
We start with the case $n\ge 2$. First we will show that $M$ is generated by the set
\BE
\label{3.3}
\left\{\frac{\sin(\pi p/m)}{\sin (\pi/m)}\right\}\cup \left\{\frac{\sin(\pi k/m)}{\sin (\pi/m)}; k=2,\LD, m', (k,p)=1\right\}.
\EE
In what follows we use the abbreviations
\BD
  \prod_{j;\, p^r}\enspace \MB{ and }\enspace \sum_{j;\,p^r}
\ED
for products and sums running over all $j\in\{1,\LD,m\}$ that are $\equiv 1\mod p^r$.

\begin{prop} 
\label{p1}
Let $m=p^n$, $n\ge 2$, be as above, $r\in\{1,\LD,n-1\}$, and $k\in\Z$, $(k,p)=1$. Then
\begin{eqnarray}
\label{3.5}
  \left|\frac{\sin(\pi p^rk/m)}{\sin(\pi/m)}\right|^{p-1}=\left(\frac{\sin(\pi p/m)}{\sin(\pi/m)}\right)^{p^r-1}\prod_{j;\,p^{n-r}}\left|\frac{\sin(\pi kj/m)}{\sin(\pi/m)}\right|^{p-1} \times \nonumber\\
  \times \prod_{j;\,p^{n-1}}\left(\frac{\sin(\pi j/m)}{\sin(\pi/m)}\right)^{1-p^r}.
\end{eqnarray}

\end{prop} 

\MN
{\em Proof.} Note that all quantities on the right hand side of formula (\ref{3.5}) are positive. From (\ref{3.1}) we see that
\BD
\label{3.6}
\left|\frac{\sin(\pi k/m)}{\sin(\pi/m)}\right|=\left|\frac{1-\zeta_m^k}{1-\zeta_m}\right|,\enspace k\in\Z.
\ED
Hence it suffices to show
\begin{eqnarray}
\label{3.7}
  \left(\frac{1-\zeta_m^{p^rk}}{1-\zeta_m}\right)^{p-1}=\left(\frac{1-\zeta_m^p}{1-\zeta_m}\right)^{p^r-1}\prod_{j;\,p^{n-r}}\left(\frac{1-\zeta_m^{kj}}{1-\zeta_m}\right)^{p-1} \times \nonumber\\
  \times \prod_{j;\,p^{n-1}}\left(\frac{1-\zeta_m^j}{1-\zeta_m}\right)^{1-p^r}.
\end{eqnarray}
We compare the denominators on both sides of (\ref{3.7}). On the left we have $(1-\zeta_m)^{p-1}$ and on the right $(1-\zeta_m)^{p^r-1}(1-\zeta_m)^{p^r(p-1)}(1-\zeta_m)^{p(1-p^r)}$.
Here we have used
\BD
 \prod_{j;\,p^{n-r}} (1-\zeta_m)=(1-\zeta_m)^{p^r}
\ED
and the respective identity for the case $r=1$. Obviously, both denominators coincide.
In the case of the numerators we use
\BE
\label{3.9}
 \prod_{j;\,p^{n-r}}(1-\zeta_m^{kj})=1-\zeta_{p^{n-r}}^k,\enspace\MB{ for } k\in\Z,\enspace (k,p)=1,
\EE
together with the special case $r=1$, $k=1$ of this formula.
Formula (\ref{3.9}) follows from the identity of polynomials
\BD
  \prod_{j;\,p^{n-r}}(X-\zeta_m^{kj})=\prod_{l=0}^{p^r-1}(X-\zeta_m^k\zeta_{p^r}^{kl})=X^{p^r}-\zeta_{p^{n-r}}^k,
\ED
which, in turn, follows from $(\zeta_m^k\zeta_{p^r}^{kl})^{p^r}=\zeta_{p^{n-r}}^k$, $l=0,\LD p^r-1$ (observe that $\zeta_{p^r}^k$ is a primitive $p^r$th root of unity).
By (\ref{3.9}), the numerator on the right hand side of (\ref{3.7}) takes the form
\BD
  (1-\zeta_{p^{n-1}})^{p^r-1}(1-\zeta_{p^{n-r}}^k)^{p-1}(1-\zeta_{p^{n-1}})^{1-p^r}=(1-\zeta_{p^{n-r}}^k)^{p-1},
\ED
which is the numerator on the left hand side.
\STOP

\MN
Next we observe that
\BD
  \left|\frac {\sin(\pi k/m)}{\sin(\pi/m)}\right|=\left|\frac {\sin(\pi k'/m)}{\sin(\pi/m)}\right|
\ED
for all numbers $k,k'\in \Z, k\equiv\pm k'\mod m$. This leads to the following: To every $k\in \Z$, $k \not\equiv 0\mod m$, there is a uniquely determined integer $k_{red}$ in $\{1,\LD,m'\}$, such that
$k_{red}\equiv \pm k\mod m$. Then
\BD
\label{3.4}
  \left|\frac {\sin(\pi k/m)}{\sin(\pi/m)}\right|=\frac {\sin(\pi k_{red}/m)}{\sin(\pi/m)}.
\ED
With this notation, Proposition 1 readily gives $U$-relations in the above sense.

\begin{theorem} 
\label{t1}
Let $m=p^n$, $n\ge 2$, be as above, $r\in\{1,\LD,n-1\}$, and $k\in\Z$, $(k,p)=1$. Then
\BE
\label{3.11}
  U_{(p^rk)_{red}}=\frac{p^r-1}{p-1} U_p+\sum_{j;\,p^{n-r}}U_{(jk)_{red}}-\frac{p^r-1}{p-1}\sum_{j;\,p^{n-1}}U_{j_{red}}.
\EE
\end{theorem} 

\MN
Observe that the set $\{(kj)_{red}; j\in\{1,\LD,m\}, j\equiv 1 \mod  p^{n-r}\}$ is not always identical with  $\{j_{red}; j\in\{1,\LD,m\}, j\equiv k\mod p^{n-r}\}$.

\MN
{\em Example.}
Let $m=27$, $r=1$, $k=2$. Then (\ref{3.11}) reduces to
\BD
 U_6=U_3+U_2+U_7-U_8-U_{10}+U_{11}.
\ED

\MN
Next we observe that the generators of $M$ given in (\ref{3.3}) are independent (as generators of a multiplicative group). Indeed,
for $2\le k\le m'$, $(k,p)=1$, the numbers $\sin(\pi k/m)/\sin(\pi/m)$ form a system of independent units in the Ring $\Z[\zeta_m]$ (see \cite[p. 144 f.]{Wa}).
On the other hand, $1-\zeta_m^j$ generates the prime ideal of $\Z[\zeta_m]$ lying above $p$ for every $j\in\Z$, $(j,p)=1$. In particular, $1-\zeta_m^j$ is associated to
$1-\zeta_m$ in the ring $\Z[\zeta_m]$. Therefore,
\BD
  1-\zeta_m^p=\prod_{j;p^{n-1}}(1-\zeta_m^j)\enspace \MB{ and }\enspace (1-\zeta_m)^p
\ED
are associated, hence also $(1-\zeta_m^p)/(1-\zeta_m)$ and $(1-\zeta_m)^{p-1}$. In particular, $(1-\zeta_m^p)/(1-\zeta_m)$ is not a unit.
From (\ref{3.1}) we obtain the said independence.

The map $\log_2:\R^+\to\R: x\mapsto \log_2(x)$ is a group isomorphism between the multiplicative group of positive real numbers and the additive group $\R$. It maps the aforesaid set of independent
generators of $M$ onto the $\Z$-linearly independent  set
\BD
  \{U_p\}\cup\{U_k; k=2,\LD, m', (k,p)=1\}.
\ED
Accordingly, these $U_k$'s are also $\Q$-linearly independent. By Theorem \ref{t1}, they form a basis of the vector space $\Q U_2+\LD+\Q U_{m'}$. Therefore, the space of $U$-relations has the
$\Q$-dimension $m'-1-\phi(m)/2$. On the other hand, Theorem \ref{t1} gives the $U$-relations
\BE
\label{3.13}
  R_{r,k}=-Y_{p^rk}+\frac{p^r-1}{p-1} Y_p+\sum_{j;\,p^{n-r}}Y_{(jk)_{red}}-\frac{p^r-1}{p-1}\sum_{j;\,p^{n-1}}Y_{j_{red}},
\EE
for $r=1,\LD, n-1$, $1\le k\le p^{n-r}/2$, $(k,p)=1$, with $k>1$ if  $r=1$. Here we use the convention $Y_1=0$ again and observe that $(p^rk)_{red}=p^rk$ for these values of $r$ and $k$.
Counting the relations $R_{r,k}$ yields a total of
$(p^{n-1}-3)/2$ relations of this kind if $p\ge 3$, and $2^{n-2}-1$ relations if $p=2$. These numbers coincide with the above $\Q$-dimension $m'-1-\phi(m)/2$.
The relations $R_{r,k}$ are $\Q$-linearly independent. Indeed, the variable $Y_{p^rk}$ occurs in $R_{r,k}$ with the coefficient $-1$,
since the numbers $(jk)_{red}$ and $j_{red}$ are not divisible by $p$ for the respective values of $k$ and $j$.
On the other hand, $Y_{p^rk}$ does not occur in any other $R_{r',k'}$ for $(r',k')\ne (r,k)$.
Altogether, we have shown

\begin{theorem} 
\label{t2}
Let $m=p^n$, $p$ a prime, $n\ge 2$.
The $U$-relations $R_{r,k}$ of {\rm (\ref{3.13})} with $r=1,\LD,n-1$, $1\le k\le p^{n-r}/2$, $(k,p)=1$, and $k>1$ if $r=1$, form a basis of the $\Q$-vector space of $U$-relations.
\end{theorem} 

\MN
Finally, we turn to the case $n=1$. So $m=p$ is a prime $\ge 5$ and $m'=(p-1)/2$. The group $M$ is generated, by its definition, by the numbers $\sin(\pi k/p)/\sin(\pi/p)$,
$k=2,\LD,(p-1)/2$. Since these numbers form a system of independent units in the Ring $\Z[\zeta_p]$ (see \cite[p. 144 f.]{Wa}), the numbers $U_k$, $k=2,\LD,(p-1)/2$, are $\Q$-linearly independent.
So there is only the trivial $U$-relation $0$.

The isomorphism $\phi$ of the foregoing section allows the computation of a basis of the $\Q$-vector space of $S$-relations in any given case (see (\ref{2.11})). In particular, we obtain the following
examples.

{\em Examples.} 1. Let $m=27$. The $U$-relations $R_{1,2}$, $R_{1,4}$ and $R_{2,1}$ forming a basis of the space of $U$-relations are transformed into the following $S$-relations:
\begin{eqnarray*}
&&  X_1+X_2-X_4-2X_5-2X_6-3X_7-3X_8-3X_9-2X_{10}-2X_{11}-2X_{12},\\
&&  X_1+2X_2+2X_3+X_4-X_5-3X_6-5X_7-6X_8-7X_9-7X_{10}-7X_{11}-6X_{12},\\
&&  3X_1+5X_2+3X_3-4X_5-8X_6-13X_7-15X_8-16X_9-14X_{10}-13X_{11}-12X_{12}.
\end{eqnarray*}
Inserting $S_1,\LD,S_{12}$ into these relations gives $0$. By Gauss elimination, we obtain
\begin{eqnarray*}
S_1&=&S_4+3S_5+S_6+S_7+S_9+S_{10}+3S_{11}+2S_{12},\\
S_2&=&-S_5+S_6+2S_7+3S_8+2S_9+S_{10}-S_{11},\\
S_3&=&-S_4+S_9+2S_{10}+3S_{11}+2S_{12}.
\end{eqnarray*}
This shows that $S_4,\LD,S_{12}$ form a $\Q$-Basis of $\Q S_1+\LD+\Q S_{12}$.

2. Let $m=32$. In this case we obtain, by the same procedure,
\begin{eqnarray*}
S_1&=&S_8+2S_9+2S_{10}+2S_{11}+4S_{12}+5S_{13}+7S_{14}+8S_{15},\\
S_2&=&2S_{10}+4S_{11}+2S_{12}+2S_{13}+S_{14},\\
S_3&=&-S_{11}+S_{12}+2S_{13}+3S_{14}+4S_{15},\\
S_4&=&S_8+2S_9,\\
S_5&=&-S_9+S_{10}+2S_{11}+S_{12},\\
S_6&=&-S_8+S_{12}+2S_{13}+S_{14},\\
S_7&=&S_{14}+2S_{15}.
\end{eqnarray*}
Again, $S_8,\LD,S_{15}$ form a $\Q$-basis of $\Q S_1+\LD +\Q S_{15}$

\MN
By a result of Baker, the $\Q$-bases in the above examples are linearly independent over the field of algebraic numbers (see \cite[Ch. 3]{Ba}). If Schanuel's conjecture is true, each of these bases even consists of
algebraically independent elements up to at most one (see \cite[Ch. 4]{Ch}).

\section*{4. The case $m=pq$}

Let $m=pq$, where $p$ and $q$ are distinct primes $\ge 3$.  As above, put $m'=\lfloor m/2\rfloor$. Let $(k,m)=1$. In this case the $U$-relations come from
\BE
\label{4.1}
  \prod_{j;p} (1-\zeta_m^{kj})=\frac{1-\zeta_p^k}{1-\zeta_p^{kq^*}}
\EE
where the product runs over all $j\in\{1,\LD,m\}$, $j\equiv 1 \mod p$, $(j,q)=1$, and $q^*$ is defined by $qq^*\equiv 1\mod p$.
Combined with the case $k=1$ of (\ref{4.1}), this gives
\BE
\label{4.3}
 \sum_{j;p}U_{(kj)_{red}}-\sum_{j;p} U_{j_{red}}=U_{(qk )_{red}}-U_{(qq^*k)_{red}}+U_{(qq^*)_{red}}-U_q.
\EE
Here the sums run over the same numbers as the product in (\ref{4.1}). Further, $k_{red}$ is defined as in the foregoing section.

In order to obtain a suitable set of $U$-relations, we denote by $\H_p$ a subset of $\Z$ with the following properties:
For every $k\in\H_p$, $(k,q)=1$; for every $j\in\{1,\LD,(p-1)/2\}$ there is a $k\in \H_p$ such that $j\equiv k \mod p$; finally,
$k\not\equiv \pm k'\mod p$ for $k,k'\in\H_p$, $k\ne k'$. For instance, if $p<q$ we may choose $\{1,\LD,(p-1)/2\}$ as a suitable set $\H_p$.
We may assume that $1\in \H_p$. Then we denote the set $\H_p\sminus\{1\}$ by $\H_p'$. Obviously, $|\H_p'|=(p-3)/2$.

For every $k\in \H_p'$, formula (\ref{4.3}) gives the $U$-relation
\BD
\label{4.5}
 R_{p,k}=\sum_{j;p}Y_{(kj)_{red}}-\sum_{j;p} Y_{j_{red}}-Y_{(qk )_{red}}+Y_{(qq^*k)_{red}}-Y_{(qq^*)_{red}}+Y_q.
\ED
Observe the convention $Y_1=0$. In the same way we obtain, for every $l\in \H_q'$, the $U$-relation
\BD
\label{4.7}
  R_{q,l}=\sum_{j;q}Y_{(lj)_{red}}-\sum_{j;q} Y_{j_{red}}-Y_{(pl )_{red}}+Y_{(pp^*l)_{red}}-Y_{(pp^*)_{red}}+Y_p,
\ED
where $pp^*\equiv 1\mod q$.
Altogether, we have found $|\H_p'|+|H_q'|=(p+q)/2-3$ $U$-relations.

The main task of this section consists in showing that these $U$-relations are $\Q$-linearly independent. To this end let $k\in \H_p$.
Then the numbers $(kj)_{red}$, $j\in\{1,\LD,m\}$, $j\equiv 1\mod p$, $(j,q)=1$, are all distinct. Otherwise, we have the congruence
$kj\equiv\pm kj'\mod m$ with $j'\in\{1,\LD,m\}$, $j'\equiv 1\mod p$, $(j',q)=1$. In particular, $j\equiv j'\mod p$. Since $(k,m)=1$,
we have $j\equiv \pm j'\mod q$ . But this requires  $j\equiv j'\mod q$, since $j\equiv -j'\mod q$ together with $j\equiv j'\mod p$
makes both $j\equiv j'\mod m$ and $j\equiv -j'\mod m$ impossible. Altogether, we see that $j\equiv j'\mod m$ and $j=j'$.

Since the said numbers $(kj)_{red}$ are distinct, we may define, for $k\in \H_p$,
\BD
 C_k=\{(kj)_{red};\enspace j\in\{1,\LD,m\},\enspace j\equiv 1\mod p,\enspace (j,q)=1\}
\ED
and write, for $k\in \H_p'$,
\BE
\label{4.9}
   R_{p,k}=\sum_{j\in C_k}Y_j-\sum_{j\in C_1} Y_j-Y_{(qk )_{red}}+Y_{(qq^*k)_{red}}-Y_{(qq^*)_{red}}+Y_q.
\EE
In the same way we may define, for $l\in \H_q$,
 \BD
 C_l'=\{(lj)_{red};\enspace j\in\{1,\LD,m\},\enspace j\equiv 1\mod q,\enspace (j,p)=1\}
\ED
and write, for $l\in \H_q'$,
\BE
\label{4.11}
   R_{q,l}=\sum_{j\in C_l'}Y_j-\sum_{j\in C_1'} Y_j-Y_{(pl )_{red}}+Y_{(pp^*l)_{red}}-Y_{(pp^*)_{red}}+Y_p.
\EE
Observe that only the sets $C_1$ and $C_1'$ contain $1$.
The following lemma plays a key role in the proof of the linear independence of our $U$-relations.

\begin{lemma} 
\label{l1}
Let $k\in\H_p$ and $l\in \H_q$.
Then $C_k\cap C_l'\ne \emptyset$.
\end{lemma} 

\MN
{\em Proof.} By the Chinese remainder theorem, there is a number $r\in\{1,\LD,m\}$ such that
\BE
\label{4.13}
  r\equiv k \mod p,\enspace r\equiv l\mod q.
\EE
If $r>m/2$, we replace $r$ by $m-r$ and have the congruences
\BD
 m-r\equiv -k \mod p,\enspace m-r\equiv -l\mod q
\ED
instead.
It suffices to consider one of these two cases. Accordingly, let $r\le m/2$ and (\ref{4.13}) hold.
We show that $r\in C_k\cap C_l'$. To this end we construct a $j\equiv 1\mod p$, $(j,q)=1$, such that $r\equiv kj\mod m$.
From (\ref{4.13}) we infer $(r,m)=1$, since $(k,p)=1$ and $(l,q)=1$.
Suppose that $r=k+sp$, $s\in\Z$. Let $t\in\Z$ be such that $s\equiv kt \mod q$ (observe $(k,q)=1$).
Then
\BD
 r=k+sp\equiv k+ktp\equiv k(1+tp)\mod q.
\ED
Let $j\in\{1,\LD,m\}$ be such that $j\equiv 1+tp\mod m$. Then $j\equiv 1\mod p$ and $(j,q)=1$, since $kj\equiv r \mod q$ and $(r,q)=1$.
Further, $r\equiv kj\mod p$ and $r\equiv kj\mod q$. Accordingly, $(kj)_{red}=r$, since $r\le m/2$.

In the same way we find a $j'\in\{1,\LD,m\}$ such that $j'\equiv 1\mod q$, $(j',p)=1$, with $r\equiv lj'\mod m$. Hence $(lj')_{red}=r$.
\STOP

\MN
Observe that for $k,k'\in\H_p$, $k\ne k'$, $C_k\cap C_{k'}=\emptyset$. Indeed, let $j\equiv 1\mod p$, $j'\equiv 1\mod p$ be such that
$(kj)_{red}=(k'j')_{red}$. Then $kj\equiv \pm k'j'$ mod $m$, which implies $k\equiv \pm k'\mod p$ and $k=k'$.

With these tools at hand we are able to prove the following theorem.

\begin{theorem} 
\label{t3}
The $U$-relations $R_{p,k}$ of {\rm (\ref{4.9})}, $k\in\H_p'$, together with  the $U$-relations $R_{q,l}$ of {\rm (\ref{4.11})}, $l\in \H_q'$,
are $\Q$-linearly independent.
\end{theorem} 

\MN
{\em Proof.}
Suppose that $p>q\ge 3$. We treat the case $q=3$ first. In this case $\H_q'=\emptyset$, so we have only to show that the $R_{p,k}$, $k\in \H_p'$, are $\Q$-linearly independent.
Since $j\in C_k$, $k\in \H_p$, implies $(j,m)=1$, but $(qk)_{red}\equiv (qq^*k)_{red}\equiv (qq^*)_{red}\equiv q\equiv 0\mod q$, it suffices to consider
the equation
\BD
 \sum_{k\in \H_p'}c_k\left(\sum_{j\in C_k} Y_j-\sum_{j\in C_1}Y_j\right)=0.
\ED
However, the sets $C_k$, $k\in \H_p$, are pairwise disjoint. Hence this equation implies $c_k=0$ for all $k\in\H_p'$.

Next suppose that $q\ge 5$, so $\H_q'\ne \emptyset$. By the above argument, it suffices to consider the equation
\BE
\label{4.15}
 \sum_{k\in \H_p'}c_k\left(\sum_{j\in C_k} Y_j-\sum_{j\in C_1}Y_j\right)+\sum_{l\in \H_q'}c'_l\left(\sum_{j\in C_l'} Y_j-\sum_{j\in C_1'}Y_j\right)=0.
\EE
Fix $k\in\H_p'$ and $l\in\H_q'$ for the time being. Let $j\in C_k\cap C_l'$. Such a $j$ exists, by Lemma \ref{l1}. Recall that $j$ is not contained in any other $C_{k'}$, nor in any other $C_{l'}'$.
Therefore, the coefficient of the variable $Y_j$ in (\ref{4.15}) is $c_k+c_l'$, which must be $0$. This holds for arbitrarily chosen numbers $k\in\H_p'$ and $l\in \H_q'$. Accordingly, there is a number $c\in \Q$
such that $c_k=c$ for all $k\in\H_p'$ and $c_l'=-c$ for all $l\in \H_q'$. Next let $k\in \H_q'$. Then there is a number $j\in C_k\cap C_1'$, by Lemma \ref{l1}. Observe that $j\ne 1$.
This implies that the coefficient of
$Y_j$ in (\ref{4.15}) equals $c_k-\sum_{l\in \H_q'}c_l'\enspace (=0)$. In other words, $(1+|\H_q'|)c=0$. Accordingly, $c=0$, and all coefficients in (\ref{4.15}) are $0$.
\STOP

\MN
{\em Remark.} In this proof we have used the left part
\BD
\sum_{j\in C_k}Y_j-\sum_{j\in C_1} Y_j
\ED
of the relation $R_{p,k}$ in order to show the linear independence. One may ask whether the simpler looking right part $-Y_{(qk )_{red}}+Y_{(qq^*k)_{red}}-Y_{(qq^*)_{red}}+Y_q$
also works. The example $p=5$, $q=11$, $k=2$, shows that this is not the case, since it gives $-Y_{22}+Y_{22}-Y_{11}+Y_{11}=0$.

\MN
Our next goal consists in showing that the space $\Q U_2+\LD +\Q U_{m'}$ has  a $\Q$-dimension $\ge m'-1-((p+q)/2-3)=\phi(m)/2+1$. To this end we consider the multiplicative
subgroup $M$ of $\Z[\zeta_m]^{\times}$ that is generated by $\sin(\pi k/m)/\sin(\pi/m)$, $k=2,\LD,m'$. The group $M$ contains the elements
\BE
\label{4.17}
 \frac{\sin(\pi k/m)}{\sin\pi/m}\left|\frac{\sin(\pi k/p)}{\sin(\pi/p)}\right|\left|\frac{\sin(\pi k/q)}{\sin(\pi/q)}\right|,\, k=2,\LD, m', (k,m)=1.
\EE
This is easy to see, for instance,
\BD
\left|\frac{\sin(\pi k/p)}{\sin(\pi/p)}\right|=\frac{\sin(\pi k_{red}/p)}{\sin(\pi/m)}\cdot \frac{\sin(\pi/m)}{\sin(\pi/p)}.
\ED
The numbers of (\ref{4.17}) are the absolute values of the Ramachanda units in $\Q(\zeta_m)$ (see \cite[p. 147]{Wa}). Since these units are independent (as elements of $\Z[\zeta_m]^{\times}$),
the numbers of (\ref{4.17}) are also independent units of $\Z[\zeta_m]$. In addition, $\sin(\pi/p)/\sin(\pi/m)$ is associated to $(1-\zeta_p)/(1-\zeta_m)$, by (\ref{3.1}). However, $1-\zeta_p$ generates
the prime ideal of $\Z[\zeta_p]$
lying above $p$, whereas $1-\zeta_m$ is a unit in $\Z[\zeta_m]$. In particular, $\sin(\pi/p)/\sin(\pi/m)$ is not a unit, so it does not lie in the group generated  by the units of (\ref{4.17}).
The same argument holds for $\sin(\pi/q)/\sin(\pi/m)$. Moreover, the set of the units of (\ref{4.17}), extended by these two numbers, is also independent. Hence the group $M$ contains
a set of $\phi(m)/2-1+2=\phi(m)/2+1$ independent elements.

As in Section 3, the group isomorphism $\log_2:\R^+\to \R$ transforms this set into a $\Q$-linearly independent subset of $\Q U_2+\LD +\Q U_{m'}$. Accordingly, this vector space has a $\Q$-dimension $\ge \phi(m)/2+1$.
This means that the space of $U$-relations has a $\Q$-dimension $\le m'-1-(\phi(m)/2+1)=(p+q)/2-3$. By Theorem \ref{t3}, we know a set of $|\H_p'|+|\H_q'|=(p+q)/2-3$ $\Q$-linearly independent $U$-relations. Thus
we have shown the following.

\begin{theorem} 
\label{t4}
In the setting of Theorem {\rm \ref{t3}},
the $U$-relations $R_{p,k}$, $k\in\H_p'$, together with $R_{q,l}$, $l\in\H_q'$,
form a basis of the $\Q$-vector space of $U$-relations.
\end{theorem} 

\MN
{\em Example.} Let $p=7$, $q=5$. Then we obtain the $U$-relations
\begin{eqnarray*}
R_{7,2}&=& Y_2+2Y_5-Y_6-Y_8+Y_9-Y_{10}+Y_{12}-Y_{13}-Y_{15}+Y_{16},\\
R_{7,3}&=& Y_3+Y_4+Y_5-Y_6-Y_8+Y_{10}+Y_{11}-Y_{13}-2Y_{15}+Y_{17},\\
R_{5,2}&=& Y_2+Y_3-Y_4-Y_6+2Y_7+Y_8-Y_9-Y_{11}+Y_{12}+Y_{13}-2Y_{14}-Y_{16}+Y_{17}.
\end{eqnarray*}
As in the foregoing section, this gives
\begin{eqnarray*}
S_1&=&-S_4+S_5+2S_6+3S_7+3S_8+5S_9+4S_{10}+2S_{11}+2S_{12}+S_{13}-S_{15},\\
S_2&=&S_4-S_5-3S_6-3S_7-S_8-2S_9+2S_{11}+3S_{12}+6S_{13}+7S_{14}+8S_{15}+6S_{16},\\
S_3&=&-S_4+S_5+3S_6+3S_7+S_8+S_9-2S_{13}-2S_{14}-3S_{15}-2S_{16}.
\end{eqnarray*}
Hence $S_4,\LD,S_{16}$ form a $\Q$-basis of $\Q S_1+\LD +\Q S_{16}$.

\MN
{\em Remark.} Most of the products in this paper have to do with norms between cyclotomic fields. For instance, the left hand side of (\ref{4.1}) can be written as
\BD
  N_{\Q(\zeta_m)/\Q(\zeta_p)}(1-\zeta_m^k).
\ED
We found, however, that the concept of norm does not simplify the present subject matter.

\section*{5. The case $m=2p$}

Let $p\ge 3$ be a prime and $m=2p$. In particular, $m'=p$.
We briefly give the facts in this case. The space $\Q U_2+\LD+\Q U_p$ has the $\Q$-basis  $U_k$, $k$ even, $2\le k\le p-1$, together with $U_p$.
Moreover, the remaining values $U_k$, $3\le k\le p-2$, $k$ odd, can be expressed by this basis in the following form:
\BE
\label{5.0}
 U_k=\begin{cases}
                  U_{p-1}+U_{2k}-U_{p-k}-U_2, \MB{ if } k<p/2;\\
                  U_{p-1}+U_{2(p-k)}-U_{p-k}-U_2, \MB{ if } k>p/2.
    \end{cases}
\EE
Accordingly, the $U$-relations $-Y_k+Y_{p-1}+Y_{2k}-Y_{p-k}-Y_2$, $k$ odd, $3\le k<p/2$, and $-Y_k+Y_{p-1}+Y_{2(p-k)}-Y_{p-k}-Y_2$, $k$ odd, $p/2<k\le p-2$, are $\Q$-linearly independent.
Since the dimension of the space of $U$-relations equals $p-1-((p-1)/2+1)=(p-3)/2$, these $U$-relations form a $\Q$-basis of this space.

As to the proof, we remark that the numbers $\sin(\pi k/p)/\sin(\pi/p)$, $2\le k \le (p-1)/2$, are independent units in $\Z[\zeta_p]=\Z[\zeta_m]$ (see the end of Section 3). In addition, the elements
$\sin(\pi/p)/\sin(\pi/m)$ and $\sin(\pi p/m)/\sin(\pi /m)$ are associated to $1-\zeta_p$ and $2$ in this ring, respectively. This follows from (\ref{3.1}), if we observe that
$1-\zeta_m$ is a unit in $\Z[\zeta_p]$. Since $1-\zeta_p$ generates the prime ideal lying above $p$, we see that the whole system
\BE
\label{5.1}
   \frac{\sin(\pi/p)}{\sin(\pi/m)}, \frac{\sin(\pi p/m)}{\sin(\pi /m)}, \frac{\sin(\pi k/p)}{\sin(\pi/p)},\enspace  2\le k \le (p-1)/2,
\EE
is independent. However,
\BD
  \frac{\sin(\pi k/p)}{\sin(\pi/p)}=\frac{\sin(\pi 2k/m)}{\sin(\pi/m)}\cdot\frac{\sin(\pi/m)}{\sin(\pi/p)}, \enspace  2\le k \le (p-1)/2.
\ED
Therefore, the system (\ref{5.1}) generates the same multiplicative group as the system
\BD
\label{5.3}
   \frac{\sin(\pi 2/m)}{\sin(\pi/m)}, \frac{\sin(\pi p/m)}{\sin(\pi /m)}, \frac{\sin(\pi 2k/m)}{\sin(\pi/m)},\enspace  2\le k \le (p-1)/2.
\ED
This shows the linear independence of the $U_k$, $2\le k\le p-1$, $k$ even, together with $U_p$.
As concerns the first identity in (\ref{5.0}), we observe that $(1-\zeta_m^k)/(1-\zeta_m)$ differs from
\BD
 \frac{1-\zeta_m^{p-1}}{1-\zeta_m^2}\cdot\frac{1-\zeta_m^{2k}}{1-\zeta_m^{p-k}},
\ED
only by a factor that is a root of unity. A similar consideration proves the second identity.

\section*{6. Short $S$-relations}

Let $m\ge 4$ be arbitrary and, as above, $m'=\lfloor m/2\rfloor$.
Let $t$ denote the dimension of the $\Q$-vector space generated by all symmetric frequencies $S_1,S_2,\LD,S_{m'-1}$. In the cases $m=27,32,35$ 
we have seen that $S_{m'-t}, S_{m'-t-1},\LD, S_{m'-1}$ form a $\Q$-basis
of this space. The author has found that this assertion holds for each $m\le 35$ (which includes cases not considered here).
So this might be true for all $m\ge 4$. One may ask whether $S_1,\LD, S_t$ is also a $\Q$-basis of this space. This, however,
is not true in all cases.

\begin{prop} 
\label{p2}
Let $m\ge 4$ be even and $m'=m/2=3j+1$, where $j$ is a natural number $\ge 2$.
Then
\BD
  -S_{j-1}+S_{2j-2} +2S_{2j-1}=0.
\ED
If, however, $m/2=3j-1$, then
\BD
  -S_{j-1}+2S_{2j-1} +S_{2j}=0.
\ED

\end{prop}

\MN
{\em Proof.} In the first case we observe that
\BD
-S_{j-1}+S_{2j-2} +2S_{2j-1}=\log_2\left(\frac{\sin(\pi j/m')^3\sin(\pi(j-1)/m)\sin(\pi(j+1)/m)}{\sin(\pi(2j+1)/m)^2\sin(\pi j/m)^2\sin(\pi(j-1)/m')}\right).
\ED
If we use the identities $\sin(\pi j/m')=2\sin(\pi j/m)\cos(\pi j/m)$ and $\sin(\pi(j-1)/m')=2\sin(\pi(j-1)/m)\cos(\pi(j-1)/m)$, the expression inside the brackets becomes
\BE
\label{6.0}
\frac{4\cos(\pi j/m)^3\sin(\pi j/m)\sin(\pi(j+1)/m)}{\sin(\pi(2j+1)/m)^2\cos(\pi(j-1)/m)}.
\EE
We have to show that this expression equals $1$.
Now the numerator of (\ref{6.0}) can be transformed by means of $4\cos(\pi j/m)^3=3\cos(\pi j/m)+\cos(\pi 3j/m)$ and $\sin(\pi j/m)\sin(\pi(j+1)/m)=\frac 12 (\cos(\pi/m)-\cos(\pi(2j+1)/m))$.
Thereby, it becomes
\begin{eqnarray}
\label{6.1}
  \frac12(3\cos(\pi j/m)\cos(\pi/m)+\cos(\pi 3j/m)\cos(\pi/m)-3\cos(\pi j/m)\cos(\pi (2j+1)/m)\nonumber \\ -\cos(\pi 3j/m)\cos(\pi(2j+1)/m)).
\end{eqnarray}
Next we apply the formula $\cos(\alpha)\cos(\beta)=\frac 12(\cos(\alpha+\beta)+\cos(\alpha-\beta))$ to each of the summands of (\ref{6.1}).
This yields
\BD
 \cos(\pi(j-1)/m)/2+ \cos(\pi(3j-1)/m)/4-\cos(\pi (5j+1)/m)/4
\ED
for the numerator of (\ref{6.0}). The same technique proves that this expression also equals the denominator of (\ref{6.0}). The second identity of Proposition \ref{p2} can be proved in the same way.
\STOP

\MN
Proposition \ref{p2} shows that, in general, not every subset of $\{S_1,S_2,\LD,S_{m'-1}\}$ of cardinality $t$ can serve as a basis of the $\Q$-vector space generated by $S_1, S_2, \LD S_{m'-1}$. In the cases
$m=14$ and $m=26$ it implies that $S_1,S_2,\LD,S_t$ is not a basis of this space. This also happens in the case $m=20$ (which we have not treated in the present paper).

\MN
{\em Remark.} Sections 3--5 suggest the following question. Is it true that
\BD
  t=\phi(m)/2-1+|\{p;\, p \MB{ a prime},\, p\DIV m\}|
\ED
for an arbitrary $m\ge 4$?


\vspace{0.5cm}
\noindent
Kurt Girstmair            \\
Institut f\"ur Mathematik \\
Universit\"at Innsbruck   \\
Technikerstr. 13/7        \\
A-6020 Innsbruck, Austria \\
Kurt.Girstmair@uibk.ac.at

\end{document}